\documentclass[12pt,leqno]{amsart}

\usepackage[utf8]{inputenc}   
\usepackage[T1]{fontenc}      
\usepackage{geometry}         

\usepackage[francais]{babel}  

\usepackage{amssymb}
\usepackage{enumerate}

\voffset-1cm \markleft{\rm Poisson} 
\newtheorem{theorem}{Theorem}
\newtheorem{definition}[theorem]{Definition}

\newtheorem{lemma}[theorem]{Lemma}

\newtheorem{proposition}[theorem]{Proposition}

\newtheorem{defi}{D\'{e}finition}

\newcommand{\K}{\mathbb K}

\newcommand{\N}{\mathbb N}

\newcommand{\pf}{\noindent{\it Proof. }}
\newcommand{\ds}{\displaystyle}

\newcommand\ant{anti-associative }

\newcommand\bu{\bullet}
\newcommand\ww{\widetilde{\wedge}}

\newcommand{\wo}{\widetilde{\otimes} }

\usepackage{graphicx}

\title{Anti-associative algebras}
\author{Elisabeth Remm}
\date{}
\address{Universit\'e de Haute-Alsace, IRIMAS UR 7499, F-68100 Mulhouse, France.}
\email{elisabeth.remm@uha.fr}

\begin{document}

\maketitle

\noindent {\bf Abstract} An anti-associative algebra is a nonassociative algebra whose multiplication satisfies the identity $a(bc)+(ab)c=0.$   Such algebras are nilpotent. We describe the free \ant algebras with a finite number of generators. Other types of nonassociative algebras, obtained either by the polarization process, such as Jacobi-Jordan algebras, or obtained by deformation quantization, are associated with this class of algebras. Following Markl-Remm's work \cite{M-R-galgalim}, we describe the operads associated with these algebra classes and in particular the cohomology complexes related to deformations.

\tableofcontents
\section{Definition and basic properties}
Let $\K$ be a field of characteristic $0$. Recall that  an algebra $(A,\mu)$ over $\K$ (often simply called an algebra) is a $\K$-vector space $A$ equipped with a bilinear multiplication  $\mu$. 
\subsection{Definition}
\begin{definition}
We say that the algebra $(A,\mu)$ is anti-associative if the product $\mu$ satisfies the quadratic identity
\begin{equation}
\label{AA}
\mu(\mu(x,y),z)+\mu(x,\mu(y,z))=0
\end{equation}
for any $x,y,z \in A$.
\end{definition}

\medskip

Such an algebra cannot be unitary.

\noindent {\bf Examples}\begin{enumerate}
  \item The only non trivial  2-dimensional \ant algebra is isomorphic to the algebra defined by
  $$e_1^2=e_2, \ e_1e_2=e_2e_1=e_2^2=0.$$
  \item In dimension $3$ we have  the following non isomorphic nontrivial anti-associative algebras $(A,\cdot)$:
\begin{enumerate}
\item $e_i \cdot e_i=0, \ e_1 \cdot e_2=-e_2 \cdot e_1=e_3,$
\item $e_1 \cdot e_1=e_2, \ e_1 \cdot e_2=-e_2 \cdot e_1=e_3,$ 
\item $e_1 \cdot e_1=e_2,\ e_1 \cdot e_3=ae_2, \ e_3 \cdot e_1=be_2 ,\ e_3 \cdot e_3=e_2, $
\item $e_1 \cdot e_1=e_2,\ e_1 \cdot e_3=ae_2, \ e_3 \cdot e_1=be_2,$
\end{enumerate}
with $a,b \in \K$ and where $\{e_1,e_2,e_3\}$ is a basis of $A$.
 \item Let us denote by $(\mathcal{F}_{AA}(X),\cdot)$ the free anti-associative algebra on one generator. It is a graded algebra
$$\mathcal{F}_{AA}(X)= \oplus_{k \geq 0} \mathcal{F}_{AA}(X)(k),$$
where $\mathcal{F}_{AA}(X)(k)$ consists of the vector space of elements of degree $k$. The anti-associativity implies that the degree $0$ component is trivial.  Let us determine each  of these components. 
\begin{enumerate}
  \item Degree $2$: We put $X \cdot X= X^2$. Thus $\mathcal{F}_{AA}(X)(2)=\K\{X^2\}.$
  \item Degree $3$: We have $X \cdot (X^2)=-X^2 \cdot X$. We put $X^3=X \cdot X^2$. Then $X^2 \cdot X =-X^3$ and $\mathcal{F}_{AA}(X)(3)=\K\{X^3\}$.
  \item Degree $4$: $$\left\{
  \begin{array}{l} X \cdot X^3 = -X^2 \cdot X^2 \\
  X \cdot X^3 = - X \cdot (X^2 \cdot X) =(X \cdot X^2)\cdot X=-(X^2 \cdot X)\cdot X=X^2\cdot X^2.
  \end{array}
  \right.
  $$
  Then $X^2 \cdot X^2=0$ and any product of degree $4$ is zero. Then
$\mathcal{F}_{AA}(X)(4)=\{0\}$
 \end{enumerate}
 We deduce that $\dim \mathcal{F}_{AA}(X) = 3.$ It corresponds to the previous example (b). 
\end{enumerate}
\subsection{Nilpotency of \ant algebras}
Recall that an algebra $A$ is called nilpotent if there exists an integer $k$ such that all products of $k$ elements of $A$ are $0$. 
The smallest integer $k$ such that   the subspace $A^k$ generated by all the products of degree $k$ is zero is called the nilindex of $A$.
To simplify the notations, we write $xy$ the multiplication of two elements of $A$ in place of $x \cdot y$. 

\begin{proposition}
An \ant algebra is nilpotent of nilindex $4$.
\end{proposition}
\pf Let $a,b,c,d$ be elements of  the anti-associative algebra $A$. We have
$$(ab)(cd)=-a ( b (c d))= a((bc)d)=-(a(bc))d=((ab)c)d=-(ab)(cd).$$
Then $(ab)(cd)=0$ and all the other products of four elements are also $0$. 

\subsection{The associative and Lie multiplication algebras}
Let $a$ be an element of the \ant algebra $A$. We denote by $L_a$ and $R_a$ the left and right translation in $A:$
$$L_a(x)=ax, \ \ R_a(x)=xa.$$
Since $A$ is a nilpotent algebra of nilindex $4$, these linear maps are nilpotent of degree $3:$
$$L_a^3=R_a^3=0.$$
Let $\mathcal{M}(A)$ be the subalgebra of $End(A)$ generated by the linear maps $L_a$ and $R_a$ for all $a \in A$. It is an associative algebra called the multiplication algebra of $A$ which is nilpotent of nilindex 3. We have the identities
$$
\left\{
\begin{array}{l}
L_aL_b=-L_{ab},\\
R_aR_b=-R_{ba},\\
R_aL_b=-L_bR_a,\\
\end{array}
\right.
$$
resulting from the anti-associativity in $A$ then every element of $\mathcal{M}(A)$ can be written as linear combination of  $L_x, R_y, L_xR_y$ 

\medskip

\noindent{\bf Example.} Let us consider the free \ant algebra on one generator. It is of dimension $3$ and its multiplication satisfies
$$e_1 e_1=e_2, \ e_1 e_2=-e_2 e_1=e_3$$
(with $ e_1=X, e_2=X^2, e_3=XX^2$). 
The associative algebra $\mathcal{M}(A)$ is  $3$-dimensional generated by $L_{e_1},R_{e_1},L_{e_2}.$  

\medskip

We can consider also the Lie subalgebra $\mathcal{L}(A)$ of $gl(A)$ generated by the left and right multiplication of $A$. It corresponds to the Lie algebra associated with the associative algebra $\mathcal{M}(A)$. Since it is generated by the Lie brackets between these nilpotent linear maps, we have
\begin{proposition}
The Lie multiplication algebra associated with the \ant algebra is a $2$-step nilpotent Lie algebra.
\end{proposition}

\medskip

\noindent{\bf Example.} In the previous example, the Lie algebra $\mathcal{L}(A)$  is isomorphic to the $3$-dimensional Heisenberg algebra. 

\subsection{Derivations and anti-derivations}

Recall that a derivation of the algebra $A$ is a linear operator which satisfies
$$f(xy)=xf(y)+f(x)y.$$
The set of derivations of $A$, $D(A)$, is a Lie algebra. A derivation $f$ is called inner if $f \in D(A) \cap \mathcal{L}(A).$ For example,  we consider for any $x \in A$ the linear map $L_x-R_x$ which belongs to $\mathcal{L}(A)$. 
In \cite{RWass}, it is proved that this linear map is a derivation of the nonassociative algebra $A$ if and only if $A$ is weakly associative. Then if $A$ is an anti-associative algebra, $L_x-R_x$ is a derivation if and only if 
 $$x(yz)+y(zx)-y(xz)=0$$
for any $x,y,z \in A$.

\begin{definition}
Let $A$ be an algebra (not necessarily anti-associative). A linear map $f \in End(A)$ is an anti-derivation if this map satisfies
$$f(xy)+xf(y)+f(x)y=0$$
for any $x,y \in A.$
\end{definition}
We denote by $\widetilde{D}(A)$ the set of anti-derivation of $A$.
\begin{proposition}
If $A$ is a nonassociative algebra (\ant or not) then$[\widetilde{D}(A),\widetilde{D}(A)] \subset D(A).$ 
\end{proposition}
\pf Let $f,g$ be two anti-derivations of $A$.
$$f\circ g (xy)=f(g(xy))=-f(xg(y)+g(x)y)=f(x)g(y)+xf(g(y))+f(g(x))y+g(x)f(y).$$
Then
$$
\begin{array}{ll}
[f,g](xy)&=f(x)g(y)+xf(g(y))+f(g(x))y+g(x)f(y)-g(x)f(y)-xg(f(y))\\
& \ \ \ \ -g(f(x))y-f(x)g(y)\\
&=x(f(g(y))-g(f(y)))+(f(g(x))-g(f(x)))y\\
&=x[f,g](y)+[f,g](x)y
\end{array}
$$
and $[f,g]$ is a derivation of $A$. 
\begin{proposition}
Let $A$ be an \ant algebra. Then for all $x \in A$, the linear map $L_x+R_x$ is an anti-derivation of $A$.
\end{proposition}
\pf We have 

$(L_x+R_x)(yz)+(L_x+R_x)(y)z+y(L_x+R_x)(z)=x(yz)+(yz)x+(xy)z+(yx)z+y(xz)+y(zx)=0$ since $A$ is anti-associative.

\noindent {\bf Remark.} Let $A$ be an algebra (not necessarily anti-associative). Then $L_x+R_x$ is an anti-derivation if and only if the multiplication of $A$ satisfies:
\begin{equation}
\label{weakly anti ass}
\widetilde{\mathcal{A}}(x,y,z)+\widetilde{\mathcal{A}}(y,x,z) +\widetilde{\mathcal{A}}(y,z,x)=0
\end{equation}

where $\widetilde{\mathcal{A}}$ is the trilinear map
$$\widetilde{\mathcal{A}}(x,y,z)=x(yz)+(xy)z.$$
An algebra whose operator $\widetilde{\mathcal{A}}$ satisfies Identity (\ref{weakly anti ass}) will be called {\it ($Id+c+\tau_{12}$)-anti-associative.}
(the terminology comes from the analogy with $v$-associative algebra where $v$ is a vector of $\K[\Sigma_3]$, the group algebra of the symmetric group of order 3, see
\cite{G.R.Nonass}).  

\medskip

\begin{definition}
Let $A$ be an anti-associative algebra.
For $x \in A$,  the anti-derivation $\widetilde{ad}_x=L_x+R_x$  will be called {\it inner} and  
$\widetilde{ \mathcal{I} } (A)$ denotes the set of inner anti-derivations of A.
\end{definition}

For any derivation $g$ of $A$ we have
$$\begin{array}{ll}
   [\widetilde{ad}_x,g] (y)  &    = \widetilde{ad}_x(g(y))-g(xy+yx)\\
      &   = xg(y)+g(y)x-g(x)y-xg(y)-yg(x)-g(y)x \\
      & =  -g(x)y-yg(x)\\
      & =\widetilde{ad}_{-g(x)}(y).
\end{array}
$$
Then $[\widetilde{ad}_x,g]=\widetilde{ad}_{-g(x)}$ and
\begin{proposition}\label{deriv-antideriv}  Let $A$ be an anti-associative  algebra. The set of inner anti-derivations $\widetilde{\mathcal{I}}(A)$  and the set of derivations  $D(A)$ satisfy 
$$[\widetilde{\mathcal{I}}(A),D(A)] \subset \widetilde{\mathcal{I}}(A).$$
\end{proposition}

\noindent {\bf Remark}. The inner anti-derivations are defined in the same way for  $Id+c+\tau_{12}$)-anti-associative algebra and Proposition {\ref{deriv-antideriv}} is also true for this type of algebras.

\medskip

Assume now that   the linear map $L_x+R_x $ belongs to $D(A)$ in the anti-associative algebra $A$. This is equivalent to write that, for all $y,z$ in $A$,
$$x(yz)+(yz)x-y(xz+zx)-(xy+yx)z=0,$$ that is,
$$x(yz)+(yz)x=0 \Leftrightarrow \widetilde{ad}_x(yz)=0$$

Then $\widetilde{ad}_x \in D(A)$ if and only if $\widetilde{ad}_x (u)=0 \ {\rm for all} \ u \in A^2.$
We deduce:
\begin{proposition} Let $A$ be an anti-associative  algebra. The set of anti-derivations $\widetilde{\mathcal{I}}(A)$  and the set of anti-derivations  $D(A)$ satisfy 
$$\widetilde{\mathcal{I}}(A) \cap D(A)=\{\widetilde{ad}_x\ {\rm such  \ that  \ for all} \  u \in A^2, \ \widetilde{ad}_x (u)=0\}.$$
\end{proposition}
In particular for any $x \in A^2$, the inner anti-derivation $\widetilde{ad}_x$ is also a derivation.  More generally,  if $f$ belongs to  $\widetilde{D}(A) \cap D(A)$ then
$$f(xy)=xf(y)+f(x)y=-xf(y)-f(x)y$$
or equivantly
$$f(x^2)=0, \ \ xf(y)+f(x)y=0$$
for any $x,y \in A$.

\medskip

For two endomorphisms $f$ and $g$ of $A$, we define the symmetric product
$$f \bu g = f \circ g + g \circ f.$$
Then
$$\begin{array}{ll}
    \widetilde{ad}_x \bu \widetilde{ad}_y  &=L_xL_y+L_xR_y+R_xL_y+R_xR_y+L_yL_x+L_yR_x+R_yL_x+R_yR_x    \\
      &  =  - L_{xy}+L_xR_y-L_yR_x-R_{yx}-L_{yx}+L_yR_x-L_xR_y-R_{xy}             \\
      & = -L_{xy+yx}-R_{xy+yx} = -\widetilde{ad}_{xy+yx}= \widetilde{ad}_{\, -(xy+yx)}.
\end{array}
$$
Likewise, if $f \in \widetilde{D}(A)$,
$$f \bu \widetilde{ad}_x=-\widetilde{ad}_{f(x)}=\widetilde{ad}_{\, -f(x)}.$$
\begin{proposition}
Let $A$ be an \ant algebra and  $f \in \widetilde{D}(A)$. Then
\begin{enumerate}
  \item $ \widetilde{ad}(x) \bu \widetilde{ad}y  =-\widetilde{ad}(xy+yx)$ for any $x,y \in A$,
  \item $f \bu \widetilde{ad}(x)=-\widetilde{ad}(f(x))$ for any $x \in A$.
\end{enumerate}
\end{proposition}
In general, the symmetric product of two \ant derivations is not an \ant derivation. In fact, if $f,g \in \widetilde{D}(A)$,
$$
\begin{array}{lll}
 (f \bu g)(xy) & =   & f(g(xy))+g(f(xy)) =-f(g(x)y+xg(y))-g(f(x)y+xf(y) )  \\
      & =& f(g(x))y+g(x)f(y)+f(x)g(y)+xf(g(y))+g(f(x))y+f(x)g(y)\\
      && +g(x)f(y)+xg(f(y))\\
      &=&2(f(x)g(y)+ g(x)f(y)) +(f(g(x))+g(f(x))y+x(f(g(y))+g(f(y)).\\   
\end{array}
$$ and $$x(f \bu g)(y)+(f\bu g)(x) y= x(f(g(y))+g(f(y)))+(f(g(x))+g(f(x))y.$$
Thus $f \bu g$ is an \ant derivation if and only if these linear maps satisfy
$$f(x)g(y)+ g(x)f(y) =-(f(g(x))+g(f(x))y-x(f(g(y))+g(f(y)).$$
\medskip
\noindent{\bf Remark: on the inner derivations of nonassociative algebras} Usually we cannot define directly the notion of inner derivation on a nonassociative algebra. If $U$ is such algebra, we consider as above the Lie multiplication algebra $\mathcal{L}(U)$ generated by the left and right multiplication of $U$. We say that a derivation $f$ of $U$ is inner if $f \in \mathcal{L}(U)$. The set
$\mathcal{L}(U) \cap D(U)$ is an ideal of $D(U)$. If  $U$ is a Lie algebra, then $\mathcal{L}(U) \cap D(U)$ corresponds to the set of adjoint operators, that is $L_x-R_x$. For an associative algebra, an element of $\mathcal{L}(U) \cap D(U)$ is of type $L_x-R_x$ as soon as $U$ is an unitary associative algebra. It is not always the case when $U$ is not unitary. For example, if we consider the $2$-dimensional associative algebra given in a basis $\{e_1,e_2\}$ by $e_1e_1=e_2$ (other non defined products are $0$), then all derivations of $L_x-R_x$ are trivial,  $D(U)$ is the set of endomorphisms $f$ satisfying $f(e_1)=ae_1+be_2, \ f(e_2)=2a e_2$, $\mathcal{L}(U)$ is the set of endomorphism $f$ satisfying $f(e_1)=be_2, \ f(e_2)=0$ and $\mathcal{L}(U) \cap D(U)$ is the set of endomorphisms $f$ satisfying $f(e_1)=be_2, \ f(e_2)=0$. Recall also that $L_x-R_x$ is a derivation of $U$ if and only if $U$ is a weakly associative algebra. 

When $A$ is an  anti-associative algebra, the brackets $[\widetilde{ad}_x,\widetilde{ad}_y]$, $[L_x-R_x,\widetilde{ad}_y]$, $[L_x-R_x
,L_y-R_y]$ are always derivations of $A$ and also \ant derivations because any product of $4$ elements of $A$ is zero. For the same reason, $L_xR_y$ and $R_xL_y$ are derivations and anti-derivations of $A$. This implies that $L_x+R_y$ is an anti-derivation as soon as $y=x+v$ with $z(vt)=0$ for any  $z,t \in A.$ In particular $L_x+R_{x+v}$ is an anti-derivation for any $v \in A^2$ and $\mathcal{L}(A)=\widetilde{\mathcal{I}}(A)$ as soon as $R_v=0$ for any $v \in A^2$.

\section{Free \ant algebras}

We denote by $\mathcal{F_{AA}}(k)$ the free \ant algebra with $k$ generators. This algebra is graded by the degree of its elements:
$$\mathcal{F_{AA}}(k)= \oplus _{1 \leq n \leq 3} \mathcal{F_{AA}}(k)^{(n)}$$
because all products of degree $4$ are null. 
Recall that the anti-associativity prevents having elements of degree $0$. 
\subsection{The free algebra $\mathcal{F_{AA}}(1)$}
We have already described this algebra in the first paragraph. Let us recall this construction briefly. Let $X$ be the generator. Then $X^2 \in  \mathcal{F_{AA}}(1)$. For the elements of degree $3$, since we have $X\cdot X^2=-X^2\cdot X$, we deduce that $\mathcal{F_{AA}}(1)^{(3)}=\K\{X\cdot X^2\}$. Then $\dim \mathcal{F_{AA}}(1)$ and a basis is given by
$$X, \ X^2, \ X\cdot X^2$$
and the multiplication
$$X\cdot X=X^2,  \  X^2\cdot X=-X\cdot X^2.$$
To use conventional conventions we denote $e_1^{(1)}=X, e_1^{(2)}=X^2,e_1^{(3)}=X\cdot X^2.$ Then we have
$$e_1^{(1)}\cdot e_1^{(1)}=e_1^{(2)}, e_1^{(1)}\cdot e_1^{(2)}=-e_1^{(2)}\cdot  e_1^{(1)}=e_1^{(3)}.$$
The Lie algebra $D(\mathcal{F_{AA}}(1))$ is $3$-dimensional constituted of matrices
$
\begin{pmatrix}
    \alpha  &  0 & 0  \\
     \beta  & 2\alpha & 0 \\
     \gamma & 0 & 3 \alpha
\end{pmatrix}
$ while $\widetilde{D}(\mathcal{F_{AA}}(1))$ is the space of matrices $
\begin{pmatrix}
    \alpha  &  0 & 0  \\
     \beta  & -2\alpha & 0 \\
     \gamma & 0 & \alpha
\end{pmatrix}
$
and an \ant derivation is a derivation when it is trivial on $(\mathcal{F_{AA}}(1))^2=\K\{e_1^{(2)},e_1^{(3)}\}.$ Any inner \ant derivation  is collinear to $\widetilde{ad}_{e_1^{(1)}}.$

\subsection{The free algebra $\mathcal{F_{AA}}(2)$}
Let $\mathcal{F_{AA}}(2) = \oplus _{1 \leq n \leq 3} \mathcal{F_{AA}}(2)^{(n)}$ the graded decomposition of $\mathcal{F_{AA}}(2)$. If we denote by $X,Y$ the generators of this free algebra, we have
\begin{enumerate}
  \item $ \mathcal{F_{AA}}(2)^{(1)}=\K\{X,Y\}$
  \item $ \mathcal{F_{AA}}(2)^{(2)}=\K\{X^2,XY,YX,Y^2\}$
  \item $ \mathcal{F_{AA}}(2)^{(3)}=\K\{X(X^2),X(XY),X(YX),X(Y^2),YX^2,Y(XY), Y(YX), YY^2\}$
\end{enumerate}
and $\dim \mathcal{F_{AA}}(2) =14.$
Moreover, we have
$$\left\{
\begin{array}{l}
  \mathcal{F_{AA}}(2)^{(1)}\cdot \mathcal{F_{AA}}(2)^{(1)} =\mathcal{F_{AA}}(2)^{(2)}   \\
   \mathcal{F_{AA}}(2)^{(1)}\cdot \mathcal{F_{AA}}(2)^{(2)} =\mathcal{F_{AA}}(2)^{(2)}\cdot \mathcal{F_{AA}}(2)^{(1)} =\mathcal{F_{AA}}(2)^{(3)}   \\\end{array}
   \right.
   $$
   all other non defined products being zero. Now let's look at the anti-derivations. Any anti-derivation of $A=\mathcal{F_{AA}}(2)$ is entirely determined by its values $f(X)$ and $f(Y)$.  In fact for any generators $uv$ of $\mathcal{F_{AA}}(2)^{(2)}$ with $u,v \in \{X,Y\}$ we have $f(uv)=-uf(v)-f(u)v.$ Any generators of $\mathcal{F_{AA}}(2)^{(3)}$  is written $u(vw)$ with $u,v,w \in \{X,Y\}$. Since $(uv)w=-u(vw)$ we have
   $$f(u(vw))=u(f(v)w), \ \ u(vf(w))-f(u)(vw)=0.$$ 
   If we take $f(X)=a_1X+a_2Y+Z_1$ and $\ f(Y)=b_1X+b_2Y+Z_2$ with $Z_i \in \mathcal{F_{AA}}(2)^{(2)}\oplus \mathcal{F_{AA}}(2)^{(3)}$, then the previous equations give $a_2=b_1=0, \ a_1=b_2.$ As a consequence for all $u \in \mathcal{F_{AA}}(2)^{(3)}$
   $$f(u)=a_1u,$$
   and for any $v \in \mathcal{F_{AA}}(2)^{(2)}$,
   $$f(v)=-2a_1v + Z, \ Z \in \mathcal{F_{AA}}(2)^{(3)}.$$
   In particular
  $$\dim \widetilde{D}(\mathcal{F_{AA}}(2))=25= 2\left(\dim  \mathcal{F_{AA}}(2)^{(2)}\oplus \mathcal{F_{AA}}(2)^{(3)}\right)+1.$$
   
   Now, concerning the inner anti-derivations, we can easily see that the linear maps $\widetilde{ad}_x$ for any $x \in \{X,Y,X^2,XY,YX,Y^2\}$ are linearly independent. We deduce that $$\dim \widetilde{I}(\mathcal{F_{AA}}(2))=6.$$
  
  \subsection{The free algebra $\mathcal{F_{AA}}(k)$}
  The free  anti-associative algebra with $k$ generators $X_1,\cdots,X_k$  admits the following grading associated with the degree of the products of generators:
  $$\mathcal{F_{AA}}(k)=\oplus _{1 \leq n \leq 3} \mathcal{F_{AA}}(k)^{(n)}$$
  and we have
  \begin{enumerate}
  \item $ \mathcal{F_{AA}}(k)^{(1)}=\K\{X_1,\cdots,X_k\},$
  \item $ \mathcal{F_{AA}}(k)^{(2)}=\K\{X_iX_j, 1 \leq i,j \leq k\}, $
  \item $ \mathcal{F_{AA}}(k)^{(3)}=\K\{X_i(X_jX_l), \ 1 \leq i,j,l \leq k\}.$
\end{enumerate}
Then
\begin{enumerate}
  \item $ \dim \mathcal{F_{AA}}(k)^{(1)}=k,$
  \item $ \dim \mathcal{F_{AA}}(k)^{(2)}=k^2, $
  \item $ \dim \mathcal{F_{AA}}(k)^{(3)}=k^3.$
\end{enumerate}
and
$$\dim \mathcal{F_{AA}}(k)=k+k^2+k^3.$$
By induction, we prove that any anti-derivation $f$ is completely defined by the vectors $f(X_i), i=1,\cdots,k$ and these vectors satisfy
$$X_i(X_jf(X_k))-f(X_i)(X_jX_k)=0$$
which implies that the restriction of $f$ to $\mathcal{F_{AA}}(k)^{(1)}$ is equal to $a\cdot Id.$ In particular
$$\dim \widetilde{D}(\mathcal{F_{AA}}(k))= 1+k\left( \dim  \mathcal{F_{AA}}(k)^{(2)}\oplus \mathcal{F_{AA}}(k)^{(3)}\right)=1+k^3+k^4.$$
Concerning the inner anti-derivations, we can easily see that the linear maps $\widetilde{ad}_x$ for any generator $x \in\mathcal{F_{AA}}(k)^{(1)}\oplus \mathcal{F_{AA}}(k)^{(2)}$ are linearly independent.   We deduce that $$\dim \widetilde{I}(\mathcal{F_{AA}}(2))=k+k^2.$$

The free algebra $\mathcal{F_{AA}}(k)$ admits a natural grading associated with the nilpotency property:
$$\mathcal{F_{AA}}(k)=\mathcal{F_{AA}}(k)^{(1)} \oplus \mathcal{F_{AA}}(k)^{(2)} \oplus \mathcal{F_{AA}}(k)^{(3)}$$
with
$$\mathcal{F_{AA}}(k)^{(1)}\cdot \mathcal{F_{AA}}(k)^{(1)}=\mathcal{F_{AA}}(k)^{(2)},\,  \mathcal{F_{AA}}(k)^{(1)}\cdot \mathcal{F_{AA}}(k)^{(2)}=\mathcal{F_{AA}}(k)^{(3)}, \mathcal{F_{AA}}(k)^{(i)}\cdot \mathcal{F_{AA}}(k)^{(3)}=0.$$
If $\{e_1,\cdots,e_k\}$ is a basis of $\mathcal{F_{AA}}(k)^{(1)}$, that is a basis of generators of the free algebra, $\{f_{i,j}=e_ie_j, 1\leq i, j \leq k\}$ a basis of $\mathcal{F_{AA}}(k)^{(2)}$, $\{g_{i,j,l}=e_if_{j,l}\}$ a basis of $\mathcal{F_{AA}}(k)^{(3)}$, then we have
$$
\left\{
\begin{array}{l}
  e_ie_j=  f_{i,j}     \\
   e_if _{j,l}=g_{i,j,l}  \\
   e_ig_{j,l,s}=g_{j,l,s}e_i=f_{i,j}f_{l,s}=f_{i,j}g_{l,s,t}=   g_{l,s,t}f_{i,j}=0\\
   f_{i,j}e_l=-e_if_{j,l}
\end{array}
\right.
$$
for all $i,j,l,s,t \in [ \! [ i,k]\! ].$

\noindent{\bf Application: Construction of \ant algebras.} From this free \ant algebra, we deduce the \ant algebras with $k$ generators:
$$A= A^{(1)} \oplus A^{(2)} \oplus A ^{(3)}$$
where $$A^{(1)}=\mathcal{F_{AA}}(k)^{(1)}, \ A^{(2)}=A^{(1)}\cdot A^{(1)}, \ A^{(3)}=A^{(1)}\cdot A^{(2)} \oplus A^{(2)}\cdot A^{(1)}.$$
This allows the construction of the \ant algebras from the data of the generators. For example, if $\dim A^{(1)}=2$, then $\dim A^{(2)} \leq 4.$ If this dimension is equal to $4$, we have $\mathcal{F_{AA}}(2)$. If this dimension is equal to $1$, then $\dim A^{(3)} \leq 2$. We complete the table of multiplication putting $f_{i,j}e_k=-e_if_{j,k}$. 

\noindent{\bf Example  $\dim  A^{(1)}=2$.} In this case $\dim A \leq 14$. If $\dim A=14$, then $A$ is isomorphic to $\mathcal{F_{AA}}(2)$. We can describe the other \ant algebras considering the dimensions of $A^{(2)}$ and $A^{(3)}$. We will not develop here a complete classification but show how to obtain it. We put $A^{(1)}=\K\{e_1,e_2\}$.
\begin{enumerate}
\item If $\dim A^{(2)}=1$, then $A^{(3)}=0.$ In fact 

1. If $A^{(2)}=\K\{f_{1,1}\}$, then 
$$e_1f_{1,1}=f_{1,1}e_1=0, \, e_2f_{1,1}=\alpha f_{1,1}e_1=0, \, f_{1,1}e_2=-\beta e_1f_{1,1}=0$$ and 
$A^{(3)}=0.$

2. If $A^{(2)}=\K\{f_{1,2}\}$, then we can assume that $e_1e_1=e_2e_2=0$ (if not we are in the previous case). Then 
$$e_1f_{1,2}=0, \ f_{1,2}e_1= 0,  \ f_{1,2}e_2= 0, \ e_2f_{1,2}= -bf_{1,2}e_2,=0$$
This implies that $A^{(3)}=0.$
We obtain the following $3$-dimensional \ant algebras:
\begin{enumerate}
  \item $e_ie_j=\alpha_{ij}f_{1,1}, \alpha_{1,1}=1, e_if_{1,1}=f_{1,1}e_i=0$.  
  \item $e_ie_i=0, e_1e_2=-e_2e_1=f_{1,2}, e_if_{1,2}=f_{1,2}e_i=0$. 
  \end{enumerate}
\item If $A^{(2)}=\K\{f_{1,1},f_{1,2}\}$ then $\dim A^{(3)}\leq 1$ and if this dimension is not $0$, then $A^{(3)}=\K\{g_{1,1,2}\}$. We obtain the following algebras
\begin{enumerate}
\item $e_1e_1=f_{1,1},\ e_1e_2=f_{1,2}, \ e_2e_1=af_{1,1}+bf_{1,2}, \  e_2e_2=cf_{1,1}+df_{1,2}.$
  \item $$\left\{
  \begin{array}{l}
  e_1e_1=f_{1,1},\, e_1e_2=f_{1,2},\, e_2e_1=af_{1,1}+bf_{1,2}, \, e_2e_2=cf_{1,1}+df_{1,2}\\
  e_1f_{1,1}=0,\, e_1f_{1,2}=g_{1,1,2},\, e_2f_{1,1}=b^2g_{1,1,2},\, e_2f_{1,2}=(a+bd)g_{1,1,2}\\
 f_{1,1}e_1=0,\, f_{1,2}e_1=-bg_{1,1,2},\, f_{1,1}e_2=-g_{1,1,2},\, f_{1,2}e_2=-dg_{1,1,2}\\
 \end{array}
 \right.
 $$
This algebra is $5$-dimensional and can be considered as the quotient of the free algebra $\mathcal{F_{AA}}(2)$ modulo the relations $YX-aX^2-bXY, \ Y^2-cX^2-dXY.$ 
\item If we assume  that $\dim A^{(2)}=2$ and that  $A$ is not isomorphic to the previous one, then such an algebra is isomorphic to the quotient of $\mathcal{F_{AA}}(2)$ modulo the relations $XY=aX^2,YX=cX^2$. In this case $A^{(3)}=0$ and we have the algebra
$$\left\{
  \begin{array}{l}
  e_1e_1=f_{1,1},\, e_1e_2=af_{1,1},\, e_2e_1=cf_{1,1}, \, e_2e_2=f_{2,2}\\
  e_if_{j,k}=0, \, f_{j,k}e_i=0.
 \end{array}
 \right.
 $$
 This algebra is 4-dimensional.
\end{enumerate}
\end{enumerate}
All the other cases can be treeted in a similar way. It will be a next work.

\section{Homology of \ant algebras}
\subsection{The standard complex of homology}  Let $A$ be a finite-dimensional  anti-associative algebra,
$T(A)=\oplus T^n(A)$
its tensorial algebra. For any $n \geq 4$, we consider the subspace $R^n(A)$ of $T^n(A)$ generated by the vectors
$$x_{i_1}\otimes \cdots \otimes x_{i_{k-1}} \otimes x_{i_k}   x_{i_{k+1}}  \otimes x_{i_{k+2}}  \otimes \cdots  \otimes  x_{i_{l-1}} \otimes x_{i_l}   x_{i_{l+1}}  \otimes x_{i_{l+2}} \otimes \cdots \otimes x_{i_{n+2}} $$
and we denote by $\widetilde{T}^n(A)$ the factor space 
$$\widetilde{T}^n(A)=\frac{T^n(A)}{R^n(A)}$$
for $n\geq 4$ and $\widetilde{T}^i(A)=T^i(A)$ for $1\leq i \leq 3$. We denote by $a\widetilde{\otimes} b$ the multiplication on $\widetilde{T}(A)= \oplus \widetilde{T}^n(A)$ deduced from the tensor product.

Let $H_*(A,A)$ be the homology group of the chain complex $(\widetilde{T}^n(A),b_n)$ where $b_n: \widetilde{T}^{n+1}(A)\rightarrow \widetilde{T}^n(A)$ is defined by

$$\begin{array}{r l }
b_n(x_{1}\widetilde{\otimes} x_2 \widetilde{\otimes}   \cdots  \widetilde{\otimes} x_{n+1} )\!  = &  x_{1} x_2\widetilde{\otimes} x_3 \widetilde{\otimes} \cdots\widetilde{\otimes} x_{n+1} + x_{1}\widetilde{\otimes} x_2 x_3 \widetilde{\otimes} x_4 \widetilde{\otimes} \cdots\widetilde{\otimes} x_{n+1}  + \cdots  \\ & + x_{1}\widetilde{\otimes}  \cdots\widetilde{\otimes}  x_{i-1 }  \widetilde{\otimes} x_i x_{i+1 }   
 \widetilde{\otimes} x_{i+2 }  \widetilde{\otimes}  \cdots \widetilde{\otimes} x_{n+1 } + \cdots \\
&  + x_{n+1}x_1 \widetilde{\otimes}  x_2 \widetilde{\otimes}  \cdots \widetilde{\otimes}  x_{n}\\
  = &\sum_{i=1}^n  x_{1}\widetilde{\otimes}  \cdots \widetilde{\otimes} x_{i-1} \widetilde{\otimes} x_i x_{i+1} \widetilde{\otimes} x_{i+2} \widetilde{\otimes} \cdots\widetilde{\otimes} x_{n+1} \\
&   +x_{n+1}x_1 \widetilde{\otimes}  x_2 \widetilde{\otimes}  \cdots \widetilde{\otimes}  x_{n}
\end{array}$$

For example 
$$\begin{array}{r l }
b_1(x_1 \widetilde{\otimes} x_2)= & x_1x_2+x_2x_1 \\
b_2 (x_1 \widetilde{\otimes} x_2 \widetilde{\otimes} x_3)= & x_1x_2 \widetilde{\otimes} x_3+x_1 \widetilde{\otimes} x_2x_3+ x_3x_1\widetilde{\otimes} x_2\\
b_3 (x_1 \widetilde{\otimes} x_2 \widetilde{\otimes} x_3\widetilde{\otimes}x_4)= & x_1x_2 \widetilde{\otimes} x_3 \widetilde{\otimes}x_4 +x_1 \widetilde{\otimes} x_2 x_3 \widetilde{\otimes}x_4+ \\
 & x_1 \widetilde{\otimes} x_2 \widetilde{\otimes}x_3x_4 + x_4x_1\widetilde{\otimes} x_2 \widetilde{\otimes} x_3\\
\end{array}$$

\begin{lemma}
For every $n\geq 2$, $b_{n-1} \circ b_n=0.$
\end{lemma}

\pf $$\begin{array}{r l }
(b_1\circ b_2) (x_{1}\widetilde{\otimes} x_2 \widetilde{\otimes}   x_{3})& = b_1(x_1x_2 \widetilde{\otimes} x_3+x_1 \widetilde{\otimes} x_2x_3+ x_3x_1\widetilde{\otimes} x_2)\\
& = (x_1x_2) x_3+x_3(x_1x_2)+ x_1( x_2x_3)+( x_2x_3)x_1 \\ 
& \qquad + (x_3x_1) x_2+x_2(x_3x_1)=0
\end{array}$$

$$\begin{array}{r l }
(b_2\circ b_3) (x_{1}\widetilde{\otimes} x_2 \widetilde{\otimes}   x_{3} \widetilde{\otimes} x_4)& = (x_1x_2)x_3 \widetilde{\otimes} x_4+(x_1x_2) \widetilde{\otimes} (x_3x_4)+ x_4(x_1x_2)\widetilde{\otimes} x_3
+ x_1(x_2 x_3) \widetilde{\otimes} x_4  \\
& +x_1\widetilde{\otimes} (x_2x_3)x_4+ (x_4 x_1)\widetilde{\otimes} (x_2x_3)+(x_1 x_2) \widetilde{\otimes}(x_3x_4)+ x_1 \widetilde{\otimes} x_2(x_3x_4)\\
& + (x_3x_4) x_1\widetilde{\otimes}x_2+(x_4x_1)x_2 \widetilde{\otimes}x_3+(x_4x_1)\widetilde{\otimes}(x_2x_3)+x_3(x_4x_1)\widetilde{\otimes}x_2\\
& =2\left[   (x_1x_2) \widetilde{\otimes} (x_3x_4)+  (x_4 x_1)\widetilde{\otimes} (x_2x_3)  \right]=0
\end{array}$$
and for the general case

\noindent $(b_{n-1}\circ b_n) (x_{1}\widetilde{\otimes}  \cdots \widetilde{\otimes} x_{n+1})$
$$\begin{array}{l}
= \ds\sum_{i=1}^n b_{n-1} (x_{1}\widetilde{\otimes}  \cdots \widetilde{\otimes} x_{i-1} \widetilde{\otimes} x_i x_{i+1} \widetilde{\otimes} x_{i+2} \widetilde{\otimes} \cdots\widetilde{\otimes} x_{n+1} )
+b_{n-1}(x_{n+1}x_1 \widetilde{\otimes}  x_2 \widetilde{\otimes}  \cdots \widetilde{\otimes}  x_{n})\\
= \ds\sum_{i=2}^{n} x_1\wo \cdots \wo x_{i-1}(x_ix_{i+1})\wo\cdots\wo x_{n+1}+  \ds\sum_{i=1}^{n-1}x_1\wo \cdots \wo (x_{i}x_{i+1})x_{i+2}\wo\cdots\wo x_{n+1} \\
 \ \ + x_{n+1}(x_1x_2)\wo x_3 \wo \cdots \wo x_n +(x_{n+1}x_1)x_2 \wo x_3 \cdots \wo x_n +(x_nx_{n+1})x_1 \wo x_2 \wo \cdots \wo x_{n-1} \\
\  \  + x_n(x_{n+1}x_1) \wo x_2 \wo \cdots \wo x_{n-1} \\
 =0
\end{array}$$

\subsection{Homology of the free \ant algebra on one generator}
Let $\mathcal{F_{AA}}(1)$  be the free \ant algebra on one generator.This algebra  admits a basis $\{e_1,e_2,e_3\}$ with
  $$e_1e_1=e_2, \ e_1e_2=-e_2e_1=e_3$$
  all the other non defined products are zero. We have
  \begin{enumerate}
  \item $\widetilde{T}^1(\mathcal{F_{AA}}(1))=\mathcal{F_{AA}}(1) $ 
    \item $\widetilde{T}^2(\mathcal{F_{AA}}(1))=\K\{e_1\wo e_2, e_2\wo e_1, e_1\wo e_3, e_3\wo e_1, e_1 \wo e_1 \}$
    
    In fact $e_2 \wo e_2 =e_1e_1 \wo e_1 e_1 =0$, $e_2 \wo e_3=e_1e_1 \wo e_1e_2=0$ and also $e_3 \wo e_2 =e_3 \wo e_3=0.$

  \item $\widetilde{T}^3(\mathcal{F_{AA}}(1))=\K\{e_1 \wo e_1 \wo e_i, i=1,2,3, e_1 \wo e_2\wo e_1, e_1\wo e_3\wo e_1,e_2 \wo e_1 \wo e_1, e_3\wo e_1\wo e_1\}$
   \item  for  $ n \geq 3$, 
  
   $\widetilde{T}^n(\mathcal{F_{AA}}(1))=
   \K\{e_{i_1}\wo \cdots \wo e_{i_n},  j\in [\! [1,n]\! ],(i_1,\cdots, \widehat{i_j} , \cdots,i_n)=(1,1,\cdots,1), i_j \in [[1,3]] \}$
\end{enumerate}
To simplify the notations, we denote 
$$\left\{ 
\begin{array}{l}
e^n_1= e_{1}\wo \cdots \wo e_{1} \  (n \ {\rm factors})\\
e^{n,k}_i= e_{1}\wo \cdots \wo e_{1} \wo e_i \wo e_1 \wo \cdots \wo e_1,  \  i=2,3,\  k=1, \cdots, n, \ k \  {\rm  is  \ the \ place \ of \ } e_i.\\
\end{array}
\right.
$$
With this notation $\widetilde{T}^n(\mathcal{F_{AA}}(1))$ is generated by the vectors $e^n_1, e^{n,k}_i$.

\medskip

We deduce in particular \begin{enumerate}
  \item ${\rm Im} (b_1)=\K\{e_2\}$,
  \item ${\rm Im} (b_2)=\K\{e_1\wo e_2+2e_2 \wo e_1,e_1 \wo e_3-e_3 \wo e_1\}=\K\{ e^{2,2}_2+2e^{2,1}_2, e^{2,2}_3-e^{2,1}_3\}$,
  \item ${\rm Ker}(b_1)=\K\{e_1 \wo e_2, e_2 \wo e_1, e_1 \wo e_3, e_3 \wo e_1\}=\K\{ e^{2,2}_2, e^{2,1}_2,e^{2,2}_3, e^{2,1}_3\}$
\end{enumerate}
and we have
$$\dim H_1(\mathcal{F_{AA}}(1),\mathcal{F_{AA}}(1))=2.$$

More generally
\begin{enumerate}
  \item $b_n(e^{n+1}_1)= 2 e^{n,1}_2+\ds \sum_{k=2}^n e^{n,k}_2$
    \item $b_n(e^{n+1,k}_2)= e^{n,k-1}_3-e^{n,k}_3, \ k \in [[ 2, \cdots , n]],$
    \item $b_n(e^{n+1,n+1}_2)=e^{n,n}_3-e^{n,1}_3,$
    \item $b_n(e^{n+1,1}_2)= 0,$
    \item $b_n(e^{n+1,k}_3)= 0, \ k\geq 1.$
    \end{enumerate}

We deduce 
$$\dim {\rm Im} (b_n)=n. $$
The vector space ${\rm Ker}(b_{n-1})$ is generated by the vectors $e^{n,k}_3$ for $k= 1, \cdots , n$, $e^{n,1}_2$ and $\sum_{k=2}^n e^{n,k}_2$. Then
$$\dim {\rm Ker}(b_{n-1})=n+2.$$
We deduce 
\begin{lemma}
$$\dim H_n(\mathcal{F_{AA}}(1),\mathcal{F_{AA}}(1))=2.$$
A basis of this space is given by the classes of the vectors $e^{n+1,1}_2$ and $e^{n+1,1}_3$.
\end{lemma}

\section{Commutative and anti-commutative \ant algebras}
 \subsection{Commutative case}
 Let $A$ be an \ant algebra. Assume that the product is commutative. Then
 $$a(bc)=-(ab)c=-c(ab)=-c(ba)=(ca)b=b(ca)=-(bc)a=-a(bc)=-(cb)a=c(ba)$$
  and $$a(bc)=0.$$
  We deduce that any commutative \ant algebra is $3$-nilpotent. It is also a commutative associative algebra.  If for any $x \in A$ we have $x^2=0$, then $xy=0$ for any $x,y \in A$ and $A$ is abelian.
  The classification of nilpotent associative algebras is established up to the dimension $5$ \cite{Kay}. We deduce from these lists the following examples:
  \begin{enumerate}
  \item Dimension $2$: $$e_1e_1=e_2.$$
  \item Dimension $3$: $$e_1e_2=e_2e_1=e_3$$
  \item Dimension $4$
  \begin{enumerate}
  \item $e_1e_1=e_2, \ e_1e_3=e_3e_1=e_4$
  \item $e_1e_1=e_2, \ e_3e_3=e_4$
  \item $e_1e_1=e_4, \ e_2e_3=e_3e_2=e_4.$
\end{enumerate}
\end{enumerate}

\medskip

The free algebra with $p$ generators $X_1, \cdots ,X_p$ admits a basis $\{X_i,1 \leq i \leq p, X_{ij}, 1 \leq i \leq j \leq p\}$ with the multiplication
$$X_iX_j=X_jX_i=X_{ij}, \ X_iX_{jk}=0, \ X_{jk}X_{rs}=0.$$
Its dimension is equal to $p(p+3)/2.$ For example, if $p=2$, this algebra is of $5$-dimensional and its structure is given by
$$e_1e_1=e_3,\ e_1e_2=e_2e_1=e_4,\ e_2e_2=e_5$$
and we have in particular $\dim H_1=6.$

\subsection{Anti-commutative case}
An \ant algebra is  anti-commutative if for any $a,b \in A$:
$$ab=-ba.$$
In particular $a^2=0$ for any $a \in A$. The first non abelian example is in dimension $3$ and given by
$$e_1e_2=-e_2e_1=e_3.$$
In this case, this algebra is also a Lie algebra isomorphic to the Heisenberg algebra. 
\begin{proposition}
An anti-commutative \ant algebra is a Lie algebra if and only if this algebra is $2$-step nilpotent.
\end{proposition}
In fact if $A$ is an anti-commutative \ant algebra. then for all $a,b,c \in A$ we have
$a(bc)=-a(cb)=-(bc)a=(cb)a=b(ca)=-b(ac)=-c(ba)=c(ab)=-(ab)c=-(ca)b=(ac)b=-(ba)c.$

\noindent We deduce
$$a(bc)+b(ca)+c(ab)=3a(bc).$$
Then $a(bc)+b(ca)+c(ab)=0$ if and only if $a(bc)=0$ for any $a,b,c \in A$.

\medskip

We deduce that any $2$-step nilpotent Lie algebra is an anti-commutative \ant algebra.
This gives a lot of examples of anti-commutative \ant algebras. 

Concerning a general classification of the finite dimensional anti-commutative \ant algebras, the previous result show that this problem is not easy to solve because   a general classification of the $2$-step nilpotent Lie algebras is still not known. In very small dimensions,  the computations are not very difficult (the problems arise as soon as the dimension is greater that $7$). We obtain that $A$ is abelian or isomorphic to
\begin{enumerate}
  \item  $e_1e_2=e_3$ in dimension $3$,
  \item $e_1e_2=e_3, e_ie_4=0$ in dimension $4$,
  \item $e_1e_2=e_3, e_1e_4=e_5$ in dimension $5$.
\end{enumerate}

We can also describe the free anti-commutative \ant algebras $\mathcal{F}_{AA_{ac}}(p)$ with $p$ generators.

1. The  free anti-commutative algebra on one generator $X$ is one-dimensional ($X^2=0$) and corresponds to the one-dimensional Lie algebra.

2. The free algebra on two generators $X,Y$ is the vector space $\K\{X,Y,XY\}$ (all products of degree $3$ are zero). It is isomorphic to the $3$-dimensional Heisenberg algebra.

3. The free algebra on $p$ generators $X_1,\cdots,X_p$ admits the basis $\{ X_i, 1 \leq i \leq p, X_iX_j, 1 \leq i < j \leq p, X_i(X_jX_k), 1 \leq i < j< k \leq p\}.$ If we denote by $e_i=X_i, e_{ij}=X_iX_j$ and $e_{i,j,k}=X_i(X_jX_k)$, the \ant multiplication is given by
$$
\left\{
\begin{array}{l}
    e_ie_j=e_{ij}, \ 1\leq i< j\leq p,   \\
    e_ie_{jk}=e_{ijk}, 1 \leq i <j< k \leq p,\\
    e_ie_{jkl}=e_{ij}e_{kl}=e_{ij}e_{kls}=e_{ijk}e_{lsr}= 0.
\end{array}
\right.
$$
If $p \geq 3$, this algebra is not a Lie algebra. For example, if $p=3$, the dimension is $7$, a basis is given by the family $\{e_1,e_2,e_3,e_{12},e_{13},e_{23},e_{123}\}$ and the multiplication:
$$
\left\{
\begin{array}{l}
   e_1e_2=e_{12},e_1e_3=e_{13},e_2e_3=e_{23},\\
 e_1e_{23}=-e_2e_{13}=e_3e_{12}=e_{123}.
\end{array}
\right.
$$
It is not $2$-step nilpotent then it is not a Lie algebra.

\section{Polarization and depolarization of \ant algebras}

Let $(A,\mu)$ be a $\K$-algebra with multiplication $\mu$. The polarization  process determines two other multiplications $\rho_\mu$ and $\psi_\mu$ on $A$ 
define by 
$$\rho_{\mu} (x,y)=\mu(x,y)+\mu(y,x), \ \ \ \psi_\mu (x,y)= \mu(x,y)-\mu(y,x)$$
for any $x,y \in A$.
The multiplication  $\rho_\mu$ is commutative and the second one, $\psi_\mu$, is anti-commutative (i.e. skew-symmetric).

The depolarization process permits to find again $\mu$ starting of a commutative multiplication $\rho$ and of an anti-commutative multiplication by
$$\mu(x,y)=\frac{1}{2}(\rho (x,y) +\psi (x,y)).$$ 
For example, if $\mu$ is associative, then $\rho_\mu$ is a Jordan multiplication and $\psi_\mu$ a Lie bracket. An example of depolarization process concerns the Poisson algebras which can be considered  as  nonassociative algebras by this process.

To simplify we will write 
$$x \bu y = \rho_\mu (x,y)=xy+yx, \ \ [x,y]=\psi_\mu (x,y)=xy-yx.$$
Then the axiom of anti-associativity is equivalent to
\begin{equation}
\label{pol}
(x\bu y)\bu z +x \bu (y\bu z )+x\bu[y,z]+z\bu[x,y]+[x,y\bu z]-[z,x\bu y]+[x,[y,z]]-[z,[x,y]]=0.
\end{equation}

\subsection{Polarization of an \ant algebra: the commutative part}
Let $A$ be an \ant algebra with multiplication  $xy$. 
If we denote by $\mathcal{AA}_\bu$ the trilinear map
$$\mathcal{AA}_\bu(x_1,x_2,x_3)=x_1\bu (x_2 \bu x_3) +  (x_1\bu x_2 )\bu x_3, $$ then 
\begin{equation} \label{equation verif anti associateur} 
\mathcal{AA}_\bu(x_1,x_2,x_3)=x_1(x_3x_2)-x_2(x_3x_1)-x_2(x_1x_3)+x_3(x_1x_2)
\end{equation}
and
we easily obtain from Equation  (\ref{equation verif anti associateur}) the following equation:
$$
\begin{array}{l}
\mathcal{AA}_\bu(x_1,x_2,x_3)+\mathcal{AA}_\bu(x_2,x_3,x_1)+\mathcal{AA}_\bu(x_3,x_1,x_2))=0.   \\
\end{array}
$$
But $\bu$ is commutative implying that
\begin{equation}
\label{bu}
 \mathcal{AA}_\bu(x_1,x_2,x_3)-\mathcal{AA}_\bu(x_3,x_2,x_1)=0.
\end{equation}
 Then considering  $\alpha_1,\alpha_2,\alpha_3,\alpha_4, \alpha_5,\alpha_6 \in \K$,  the identity $\sum_{\sigma_i \in \Sigma_3} \alpha_i \mathcal{AA}_\bu(x_{\sigma_i(1)},x_{\sigma_i(2)},x_{\sigma_i(3)})=0$, as $\alpha_5=\alpha_1+\alpha_3-\alpha_4$ and $\alpha_6=\alpha_1-\alpha_2+\alpha_3$, is equivalent to the system 
$$
\begin{array}{l}
  \alpha_1(\mathcal{AA}_\bu(x_1,x_2,x_3)+\mathcal{AA}_\bu(x_2,x_3,x_1)+\mathcal{AA}_\bu(x_3,x_1,x_2))\\+\alpha_2( \mathcal{AA}_\bu(x_2,x_1,x_3)-\mathcal{AA}_\bu(x_3,x_1,x_2))\\
\alpha_3(  \mathcal{AA}_\bu(x_3,x_2,x_1)+\mathcal{AA}_\bu(x_2,x_3,x_1)
+\mathcal{AA}_\bu(x_3,x_1,x_2))\\+\alpha_4(\mathcal{AA}_\bu(x_1,x_3,x_2)-\mathcal{AA}_\bu(x_2,x_3,x_1))=0
\end{array}
$$
which is equivalent to the system
$$
\left\{
\begin{array}{l}
\mathcal{AA}_\bu(x_1,x_2,x_3)+\mathcal{AA}_\bu(x_2,x_3,x_1)+\mathcal{AA}_\bu(x_3,x_1,x_2))=0,    \\
   \mathcal{AA}_\bu(x_1,x_2,x_3)-\mathcal{AA}_\bu(x_3,x_2,x_1)=0.
\end{array}
\right.
$$
which is clearly satisfied. 
Since $\bu$ is commutative, this system is equivalent to
\begin{equation}
\label{bu}
x_1\bu(x_2\bu x_3)+x_2\bu (x_3\bu x_1)+x_3\bu (x_1\bu x_2)=0.
\end{equation}
\begin{proposition}
Let $(A,\mu)$ be an \ant algebra. Then the commutative multiplication $\bu$ associated with $\mu$ satisfies the Jacobi identity
$$x_1\bu(x_2\bu x_3)+x_2\bu (x_3\bu x_1)+x_3\bu (x_1\bu x_2)=0.$$
Such algebras are called Jacobi-Jordan algebras (\cite{Bu,Ka}). 
\end{proposition}

\noindent{\bf Examples on Jacobi-Jordan algebras  \cite{Bu}}:
\begin{enumerate}
  \item Dimension $2$: $e_1e_1=e_2$
  \item Dimension $3$: $e_1e_1=e_2,\ e_3e_3=e_2$
  \item Dimension $4$:\begin{enumerate}
  \item $e_1e_1=e_2,\ e_1e_3=e_4$
  \item $e_1e_1=e_2,\ e_3e_4=e_2$
\end{enumerate}
\end{enumerate}

An algebraic study of these algebras can be found in \cite{Bu,Ka}.  A Jacobi-Jordan algebra is a nilalgebra of nilindex $3$. In fact for any $x \in A$ we have $x^3=0$. More precisely, in \cite{Bu}, one proves that the class of Jacobi-Jordan algebras and the class of commutative nilalgebras of nilindex at most three coincides.  This implies, from an Albert result, that
\begin{proposition}
Any finite dimensional Jordan-Jacobi algebra is nilpotent.
\end{proposition}

\begin{definition}
We will said that an algebra $(A,\mu)$ is Jacobi-Jordan admissible is the commutative algebra $(A,\bu)$ define by $x \bu y=\mu(x,y)+\mu (y,x)$ is a Jordan-Jacobi algebra.
\end{definition}
A necessary and sufficient condition for a multiplication $\mu$ to be Jacobi-Jordan algebra is
$$\sum_{\sigma \in \Sigma_3}\mathcal{AA}_\mu(x_{\sigma(1)},x_{\sigma(2)},x_{\sigma(3)})=0$$
where $\Sigma_3$ is the symmetric group of degree $3$. 

In particular, \ant algebras are Jacobi-Jordan admissible. A Jacobi-Jordan algebra which comes from an \ant algebras will be called a {\bf special Jacobi-Jordan algebra}. In the following, we will give some examples of non special Jacobi-Jordan algebras.
\begin{proposition}
If a Jacobi-Jordan algebra  is special,  then  its nilindex is smaller or equal to $4$.
\end{proposition}
\noindent {\bf Example:  the free Jacobi-Jordan algebra $\mathcal{F}_{JJ}(p)$.} Let us consider the free Jacobi-Jordan algebra $\mathcal{F}_{JJ}(p)$ on $p$ generators. For $p=1$, we have 
 $$\mathcal{F}_{JJ}(1)=\K\{X\} \oplus \K\{X^2\}$$ and the multiplication corresponds to
 $$e_1e_1=e_2.$$
 This algebra is special because it is a commutative \ant algebra and any commutative \ant algebra is a special Jordan-Jacobi algebra. 
 
\noindent  For $p=2$, we have
 $$\mathcal{F}_{JJ}(2)=\oplus_{k \geq 1} \mathcal{F}_{JJ}(2)^{(k)}$$
 where $\mathcal{F}_{JJ}(2)^{(k)}$ is the vector space of terms of degree $k$ and
 $$\mathcal{F}_{JJ}(2)^{(k)}\cdot \mathcal{F}_{JJ}(2)^{(k')}=\mathcal{F}_{JJ}(2)^{(k+k').}$$
 We have
  \begin{itemize}
  \item $\mathcal{F}_{JJ}(2)^{(1)}= \K\{X,Y\},$
  \item  $\mathcal{F}_{JJ}(2)^{(2)}= \K\{X^2,XY,Y^2\},$
   \item  $\mathcal{F}_{JJ}(2)^{(3)}= \K\{X(XY),Y(XY)\}.$
  \item  $\mathcal{F}_{JJ}(2)^{(4)}= \K\{X(XY^2)\}.$ 
 
  In fact for the terms in degree $4$ we have to compute the terms of the form $a(b(cd))$ and $ (ab)(cd)$ with $a,b,c,d =X$ or $Y$.
  We have $$
  \left\{
  \begin{array}{l}
  (ab)(cd)+c(d(ab))+d((ab)c)=0, \\
  (cd)(ab)+a(b(cd))+b((cd)a)=0.
  \end{array}
  \right.
  $$
  implying
  $$c(d(ab))+d((ab)c)=a(b(cd))+b((cd)a).$$
This identity implies when $c=d=X,a=b=Y$:
  $$X(XY^2)=Y(YX^2)$$
and when $c=X,a=b=d=Y$:
 $$Y(Y^2X)=2Y(Y(XY)).$$
 But $Y(XY^2)+Y(Y(YX))+Y(Y(XY))=0$ then
 $Y(XY^2)+2Y(Y(XY))=Y(XY^2)+Y(Y^2X)=0$ implying that 
 $Y(XY^2)=0.$
  We deduce
 $$Y(XY^2)=Y(Y(XY))=0$$
 and also
 $$X(X^2Y)=X(X(XY))=0.$$
 In summary:
 $$X(X^2Y)=X(X(XY))=X^2(XY)=Y(XY^2)=Y(Y(XY))=Y^2(XY)=0.$$
 For the other terms we have
 $$\left\{
 \begin{array}{ll}
 (XY)(XY)=-X(Y(XY))-Y(X(XY)), \\
 X^2Y^2=-2X(XY^2)=-2Y(YX^2),\\
 X(Y(XY))+X(XY^2)+X(Y(XY))=0,\\
 Y(X(XY))+Y(X(XY))+Y(YX^2)=0.
 \end{array}
 \right.$$
 This shows that  
 
 $(XY)(XY)=X(XY^2))=Y(YX^2))=-2X(Y(XY))=-2Y(X(XY))=-\frac{1}{2}X^2Y^2$ and all the terms of degree $4$ are zero or equal to $\alpha X(XY^2).$ We remark also that all the terms of partial degree $3$ (that is the maximum of the degree of $X$ and the degree of $Y$) are zero.
\end{itemize}
We deduce
$$\dim \mathcal{F}_{JJ}(2)^{(4)}=1.$$
\begin{proposition}
The free Jordan-Jacobi algebra with two generators is not special.
\end{proposition}

\medskip

Using the previous remark which precise that any term of partiel degree greater or equal to $3$ is zero, we deduce that $\mathcal{F}_{JJ}(2)^{(5)}=\{0\}$ and more generaly $\mathcal{F}_{JJ}(p)^{(k)}=\{0\}$ as soon as $k \geq 2p+1$.

\subsection{Polarization of an \ant algebra: the anti-commutative part}
Let $A$ be an \ant algebra with multiplication  $xy$ and 
$\bu$ and $[,]$ the associated commutative and anti-commutative multiplications. These two multiplications satisfy the relation (\ref{pol}) which implies 
\begin{equation}
\label{pol1}
-y\bu(z\bu x)+x\bu[y,z]+z\bu[x,y]+[x,y\bu z]-[z,x\bu y]+[x,[y,z]]-[z,[x,y]]=0.
\end{equation}

 This equation is solved as soon as $(A,\bu)$ is nilpotent with an nilindex equal to $3$ and  $[x,y]=0$ for any $x,y \in A$. Then
\begin{proposition}
Any Jacobi-Jordan algebra with a nilindex equal to $3$ is special.
\end{proposition}
We can also examine some consequences of Equation (\ref{pol}). 
Using the natural action of the symmetric group $\Sigma_3$ on the triple $(x,y,z)$, we obtain
\begin{proposition}
The commutative and anti-commutative multiplication associated with the \ant multiplication satisfy the identities
\begin{enumerate}
  \item $[x,y \bu z]+[y, z \bu x] + [z, x \bu y]=0,$
  \item $[x,[y,z]]-2[y,[z,x]]+[z,[x,y]]=x\bu (y \bu z)-z \bu (x \bu y),$
  \item $ x \bu [y,z]+[y,z \bu x]+ [x \bu y,z]=0$
\end{enumerate}
for any $x,y,z \in A$.
\end{proposition}
\noindent{\bf Remarks.}
\begin{enumerate}
  \item The identity $ x \bu [y,z]=-[y,z \bu x]- [x \bu y,z]$ implies 
  $$ x \bu [y,z] + y \bu [z,x]+z \bu [x,y]=0.$$
  Since $[x,y \bu z]+[y, z \bu x] + [z, x \bu y]=0$ and 
  $x(yz)=x \bu (y \bu z)+[x, y \bu z]+x \bu [y,z] +[x,[y,z]]$, we have
  $$x(yz)+y(zx)+z(xy)= [x,[y,z]]+[y,[z,x]]+[z,[x,y]]$$
  and the algebra $(A,[,])$ is a Lie algebra if and only if the \ant algebra $A$ is a non-commutative Jacobi-Jordan algebra. 
  \item If we denote by $f_x$ the linear map $f_x(y)=x \bu y$, then the second relation  is written
$$f_x[y,z]+[y,f_x(z)]+[f_x(y),z]=0$$
that is $f_x$ is an anti-derivation of the  algebra $(A,[,])$. If $\widetilde{\mathcal{D}er}(A,[,])$ denotes the subvector space of $End (A)$ whose elements are the anti-derivation of $(A,[,]$, then 
$$[ \widetilde{\mathcal{D}er}(A,[,]),\widetilde{\mathcal{D}er}(A,[,])] \subset \mathcal{D}er (A,[,]).$$
 If we call inner anti-derivation the anti-derivations of type $f_x$ for $x \in A$, then
$$ f_{x \bu y}=-f_x \circ f_y -f_y \circ f_x$$
and the subspace of inner derivation is a subalgebra of the Jordan algebra $End(A)$.
\end{enumerate}
\begin{proposition}
Let $A$ be a Jacobi-Jordan algebra and $f$ an anti-derivation of $A$. We can provide on the one dimensional extension $A'=A \oplus \K\{X\}$ of $A$ a Jordan-Jacobi structure $\mu$ given by
$$
\left\{
\begin{array}{l}
  \mu (x,y)=xy, \ \forall x,y \in A    \\
   \mu (x,X)=f(x)  \ \forall x \in A, \\
   \mu(X,X)=0. 
\end{array}
\right.
$$
\end{proposition}
This proposition is a consequence of Corollary 4.6  of \cite{Agore}.

\subsection{Deformation quantization of anti-Poisson algebras}
\begin{definition}
An anti-Poisson algebra is a triple $(A,\psi,\rho)$ where $A$ is a $\K$-vector space, $\psi$ an anti-commutative (or skew-symmetric) multiplication on $A$, $\rho$ a Jacobi-Jordan multiplication on $A$ such that these two bilinear maps satisfy the anti-Leibniz identity:
$$\mathcal{L}_g(\psi,\rho)(x,y,z)=\rho(\psi(x,y),z)+\psi(x,\rho(y,z))+\psi(\rho(x,z),y)=0.$$
for any $x,y,z \in A$.
\end{definition}

Recall that an antiderivation of an algebra $(A,\mu)$ is a linear map $f$ such that
$$f(\mu(x,y))+\mu(x,f(y))+\mu(f(x),y)=0$$
for any $x,y$ in $A$. The anti-Leibniz identity can be interpreted saying that for any $z \in A$, the linear map $x \rightarrow \rho(x,z)$ is an antiderivation of the algebra $(A,\psi).$

\medskip

A formal deformation of an \ant algebra $(A,\mu_0)$ is given by a family
$$\varphi_i: A \otimes A \rightarrow A, \ \ i \in \N$$
satisfying $\varphi_0=\mu_0$ and
$$(D_k): \ \ \ \ \ \sum_{i+j=k, i,j\geq 0}\varphi_i(\varphi_j(a,b),c)
+\sum_{i+j=k, i,j\geq 0} \varphi_i(a,\varphi_j(b,c))=0$$
for any $a,b,c \in A$ and for each $k \geq 1$.
Now expand $\mu(a, b)$, for $a, b \in A$, into the power series
$$\mu(a, b) = \mu_0(a, b) + t\varphi_1(a, b) + t^2\varphi_2(a,b)+ \cdots$$
for some $k$-bilinear functions $\varphi_i: A \otimes A \rightarrow A$. Then  $\mu$ is \ant if and only if $(D_k)$ are satisfied for each $k \geq 1$.

In we denote as before $\mu_0(a,b)=ab$ the multiplication of the \ant algebra $A$, then
$$(D_1): \ \ \ \ \  \varphi_1(a,bc)+a(\varphi_1(b,c))+\varphi_1(ab,c)+(\varphi_1(a,b))c=0,$$
$$(D_2): \ \ \ \ \  \varphi_1(a,\varphi_1(b,c))+\varphi_1(\varphi_1(a, b),c))+\varphi_2(a,bc)+a(\varphi_2(b,c))+\varphi_2(ab,c)+(\varphi_2(a,b))c=0.$$
To make a link with the classical cohomological approach of deformation theories, we put
$$Z^2_{\mathcal{AA}}(A,A)=\{\varphi : A \otimes A \rightarrow A, \delta^2(\varphi)(a,b,c)= \varphi(a,bc)+a(\varphi(b,c))+\varphi(ab,c)+(\varphi(a,b))c=0.\}$$
Then $(D_1)$ is satisfied as soon as $\varphi_1 \in Z^2_{\mathcal{AA}}(A,A)$. For example, if $f$ is a linear endomorphism of $A$, the bilinear map $\varphi_f(a,b)=f(ab)-f(a)b-af(b)$ is in $Z^2_{\mathcal{AA}}(A,A)$.

Assume now that $\mu_0$ is anti-commutative: $ba=-ab$ and consider the condition $D_2$ that we write 
$$\varphi_1(x_1,\varphi_1(x_2,x_3))+\varphi_1(\varphi_1(x_1, x_2),x_3))+\delta^2(\varphi_2)(x_1,x_2,x_3)=0.$$
Let us determine the constants $a_\sigma$  such that

$$\sum_{\sigma \in \Sigma_3}a_{\sigma} \delta^2(\varphi_2)(x_{\sigma(1)},x_{\sigma(2)},x_{\sigma(3)})=0.$$
If we develop this equation we obtain the non trivial solution $a_{\sigma}=1$ for every $\sigma \in\Sigma_3$.  We deduce from $D_2$ that
$$\sum_{\sigma \in \Sigma_3} \varphi_1(x_{\sigma(1)},\varphi_1(x_{\sigma(2)},x_{\sigma(3)}))+\varphi_1(\varphi_1(x_{\sigma(1)}, x_{\sigma(2)}),x_{\sigma(3)})=0$$
and $\varphi_1$ is an \ant admissible multiplication. We deduce that the commutative multiplication $\rho_1$ associated with $\varphi_1$ is an \ant multiplication. Now, from $(D_1)$ we have 
$$\sum_{\sigma \in \Sigma_3}a_{\sigma} \delta^2(\varphi_1)(x_{\sigma(1)},x_{\sigma(2)},x_{\sigma(3)})=0$$
for any $a_{\sigma}$. In particular, if
$$a_{Id}=a_{\tau(12)}=a_{c},\ a_{\tau(13)}=a_{\tau(23)}=a_{c^2}$$
where $\tau(ij)$ is the transposition between $i$ and $j$ and $c$ the cycle $(123)$, then we have
$$\varphi_1(x_1,x_2x_3)+x_2\varphi_1(x_3,x_1)-x_3\varphi_1(x_1,x_2)+x_2\varphi_1(x_1,x_3)-x_3\varphi_1(x_2,x_1)+\varphi_1(x_2x_3,x_1)=0$$
that is
$$\rho_1(x_1,x_2x_3)+x_2\rho_1(x_1,x_3)-x_3\rho_1(x_1,x_2)=0$$
or since $\mu_0$ is anti-commutative,
$$\rho_1(x_1,x_2x_3)+x_2\rho_1(x_1,x_3)+\rho_1(x_1,x_2)x_3=0.$$
We shall say that the multiplication $\rho_1$ and $\mu_0$ are connected by the anti-Leibniz relation. 
\begin{theorem}
Let $(A,\mu_0)$ be an anti-symmetric \ant algebra and $\mu_t= \mu_0+ \sum t^k \varphi_k$ an \ant formal deformation of $\mu_0$, that is $\mu_t$ is an \ant multiplication. Then $(A,\mu_0,\rho_1)$ is an anti-Poisson algebra, that is
\begin{enumerate}
\item $(A,\rho_1)$ is a Jordan-Jacobi algebra,
\item The products $\mu_0$ and $\rho_1$ are connected by the anti-Leibniz identity:
$$\mathcal{L}_g(\mu_0,\rho_1)(x,y,z)=\rho_1(xy,z)+x\rho_1(y,z)+\rho_1(x,z)y=0.$$
\end{enumerate}
We will say that $\mu_t$ is a deformation quantization of the anti-Poisson algebra $(A,\mu_0,\rho_1)$.
\end{theorem}

\section{The operad $\mathcal{AA}ss$}
\subsection{The operad $\mathcal{AA}ss$}
Recall that deformations of an associative algebra  are controlled by
the Hochschild cohomology $H^*(A,A)$. This means that $H^1(A,A)$ classifies infinitesimal deformations and
$H^2(A,A)$ contains the obstructions for their extensions. In this case this cohomology coincides with the operadic cohomology because the corresponding quadratic operad $\mathcal{A}ss$ is a Koszul operad. In this section we recall some results given in \cite{M-R-galgalim} and we refer to this paper for the all the proofs. 
Let $\mathcal{AA}ss=\oplus_{n \geq 1} \mathcal{AA}ss(n)$ the quadratic operad associated with  \ant algebras. We have\begin{enumerate}
  \item $\mathcal{AA}ss(1)=\K$
  \item $\mathcal{AA}ss(2)=\K\{x_1x_2,x_2x_1\}$
  \item $\mathcal{AA}ss(3)=\K\{x_1(x_2x_3), x_2(x_1x_3),x_3(x_2x_1),x_1(x_3x_2),x_2(x_3x_1),x_3(x_1x_2)\}$
  \item $\mathcal{AA}ss(n)=\{0\}$ for $n \geq 4$.
\end{enumerate}
The Poincar\'e or generating series of $\mathcal{AA}ss$ is the series
$$g_{\mathcal{AA}ss(t)}=\sum \frac{(-1)^n}{n!}\dim (\mathcal{AA}ss (n))t^n=-t+t^2-t^3.$$
Since $\mathcal{AA}ss(1)=\K$, $\mathcal{AA}ss$ admits a minimal model, unique up to isomorphism. The generating series $g_M$ associated with  this minimal model is the formal inverse of $g_{\mathcal{AA}ss(t)}$ taken with the opposite sign
$$g_{\mathcal{AA}ss}(-g_M(t))=t.$$
We deduce
$$g_M(t)=-t+t^2-t^3+ 4t^5-14t^6+ 30t^7 -33t^8 -55t^9+ \cdots$$
and this series cannot be the Poincare Series of the dual operad. In fact, if $\mathcal{R}$ is the ideal of relations which defines 
$\mathcal{AA}ss$ and which is generated as a $\Sigma_3$-module by the relation $(x_1(x_2x_3)+(x_1x_2)x_3)$, the ideal $\mathcal{R}^\bot$ which generates the dual operad $\mathcal{AA}ss^{!}$ is the orthogonal of $\mathcal{R}$ for the classical inner product
\begin{eqnarray}
\label{pairing}
\left\{
\begin{array}{l}
<(x_i \cdot x_j)\cdot x_k,(x_{i'} \cdot x_{j'})\cdot x_{k'}>=0, \ {\rm if} \ (i,j,k) \neq (i',j',k'), \\
<(x_i \cdot x_j)\cdot x_k,(x_i \cdot x_j)\cdot x_k>={\varepsilon(\sigma)},  
\qquad   {\rm with} \ \sigma =
\left(
\begin{array}{lll}
1 & 2 & 3\\
i &j &k 
\end{array}
\right)
\\
<x_i \cdot (x_j\cdot x_k),x_{i'} \cdot (x_{j'}\cdot x_{k'})>=0, \ {\rm if} \ ( i,j,k) \neq ( i',j',k'), \\
<x_i \cdot (x_j\cdot x_k),x_i \cdot (x_j\cdot x_k)>=-{\varepsilon(\sigma)} 
\qquad {\rm with} \ \sigma =
\left(
\begin{array}{lll}
1 & 2 & 3\\
i &j &k 
\end{array}
\right)
, \\
<(x_i \cdot x_j)\cdot x_k,x_{i'} \cdot (x_{j'}\cdot x_{k'})>=0,
\end{array}
\right.
\end{eqnarray}
where $\varepsilon(\sigma)$ is the signature of $\sigma$. Then $\mathcal{R}^\bot$ is also generated as a $\Sigma_3$-module by the relation $(x_1(x_2x_3)+(x_1x_2)x_3)$ and $\mathcal{AA}ss$  is selfdual, it was also the case for the operad $\mathcal{A}ss$ of associative algebra. We deduce that the generating series of $\mathcal{AA}ss^{!}$ is $-t+t^2-t^3$. We deduce
\begin{theorem}\cite{M-R-galgalim} . The quadratic operad $\mathcal{AA}ss$ is not a Koszul operad.
\end{theorem}
Concerning the problem of deformation of anti-associative algebras, the `standard' cohomology  of an
anti-associative algebra $A$ with coefficients in itself is described in \cite{M-R-galgalim}  and compared
to the relevant part of the deformation cohomology  based on
the minimal model of the anti-associative operad $\mathcal{AA}ss$. Since
$\mathcal{AA}ss$ is not Koszul, these two cohomologies
differ. The standard cohomology $H^*_{st}(A,A)$ is the cohomology of the complex
$$
C^1(A,A) \stackrel{\delta^1_{AA}}{\longrightarrow} C^2(A,A)
\stackrel{\delta^2_{AA}}{\longrightarrow}
C^3(A,A) \stackrel{\delta^3_{AA}}{\longrightarrow}
0 \stackrel{0}{\longrightarrow} 0 \stackrel{0}{\longrightarrow} \cdots
$$
in which $C^p(A,A) := Hom(A^{\otimes^n},A)$ for $p = 1,2,3$, and all higher
$C^p$'s are trivial. The two nontrivial pieces of the differential are
basically the Hochschild differentials with ``wrong'' signs of some
terms:
\begin{align*}
\delta^1_{AA}(f)(x,y)& := x f(y) - f(xy) + f(x)y,
\mbox { and}
\\
\delta^2_{AA}(\varphi)(x,y,z) &:= x \varphi(y,z)+ \varphi(xy,z)+ \varphi(x,yf) +\varphi(x,y)z,
\end{align*}
for $f\in Hom(V,V)$, $\varphi \in Hom(A^{\otimes^2},V)$ and $x,y,z \in
V$. One sees, in particular, that 
$$H^*_{st}(A,A)^p = 0$$ for $p\geq 4$.

\medskip

\subsection{The minimal model}
In case of \ant algebras, the deformation cohomology is described by the cohomology defined by the minimal model.This minimal model is an homology isomorphism $(\Gamma(E),\partial) \rightarrow \mathcal{AA}ss$ from the free operad $\Gamma(E)$  equipped to a differential and $\mathcal{AA}ss$. The Poincaré series, studied above as the inverse of the Poincaré series of   $\mathcal{AA}ss$, that is 
$$g_M(t)=-t+t^2-t^3+ 4t^5-14t^6+ 30t^7 -33t^8 -55t^9+ \cdots$$
corresponds to the series associated with the minimal model. 
The  deformation cohomology of anti-associative algebras, based of the study of a minimal model is also studied in \cite{M-R-galgalim}. We summarize the results:
we consider the complex 
\[
C^1_{AAss}(A,A) \stackrel{\delta^1}{\longrightarrow} C^2_{AAss}(A,A) 
\stackrel{\delta^2}{\longrightarrow}
C^3_{AAss}(A,A)  \stackrel{\delta^3}{\longrightarrow}
C^4_{AAss}(A,A)  \stackrel{\delta^4}{\longrightarrow} \cdots
\]

-- $C^1_{AAss}(A,A) = Hom(A,A)$

-- $C^2_{AAss}(A,A) = Hom(A^{\otimes^2},A)$

-- $C^3_{AAss}(A,A) = Hom(A^{\otimes^3},A)$, and

-- $C^4_{AAss}(A,A) = Hom(A^{\otimes^5},A) \oplus Hom(A^{\otimes^5},A)  \oplus
 Hom(A^{\otimes^5},A) \oplus Hom(A^{\otimes^5},A) $.

\noindent
Observe that $C^p_{AAss}(A,A)= C^p_{AA} p$ for $p = 1,2,3$, while 
$C^4_{AAss}(A,A)$ consists of $5$-linear maps. The
differential $\delta^p$ agrees with $\delta_{AA}^p$ for  $p = 1,2$ while,
for $g \in C^3_{AAss}(A,A)$, one has
\[
\delta^3 (g) = (\delta_1^3(g),\delta_2^3(g),\delta_3^3(g),\delta_4^3(g)),
\]
where
\begin{align*}
\delta_1^3(g)(x,y,z,t,u) &:= 
xg(y,z,tu) - g(x,y,z(tu)) + (xy)g(z,t,u) - g(xy,zt,u) 
\\
&\
+ g(xy,z,t)u - g((xy)z,t,u) + g(x,y,z)(tu) - g(x,yz,tu),
\\
\delta_2^3(g)(x,y,z,t,u) &:= 
g((xy)z,t,u) - g(xy,z,t)u + g(x,y,zt)u - g(x,y(zt),u)
\\
&\
+ xg(y,zt,u) - g(x,y,(zt)u) + (xy)g(z,t,u) - g(xy,z,tu),
\\
\delta_3^3(g)(x,y,z,t,u) &:= 
g(x,yz,tu) - xg(yz,t,u) + g(x,(yz)t,u) - x(g(y,z,t)u)
\\
&\
+ g(x,y,zt)u - g(xy,z,t)u + (g(x,y,z)t)u - g(x(yz),t,u), \mbox { and}
\\
\delta_4^3(g)(x,y,z,t,u) &:= 
g(xy,zt,u) - g(x,y,(zt)u) + x g(y,zt,u) - g(x,y(zt),u)
\\
&\
+ (xg(y,z,t))u - g(x,yz,t)u + (g(x,y,z)t)u - g(xy,z,t)u,
\end{align*}
for $x,y,z,t,u \in V$.  

\noindent{\bf Case of Jacobi-Jordan algebra.}  Let $\mathcal{J}a\mathcal{J}o=\oplus_{n \geq 2} \mathcal{J}a\mathcal{J}o (n)$ be the quadratic operad associated with the Jacobi-Jordan multiplication. It is clear that
\begin{enumerate}
  \item $\dim \mathcal{J}a\mathcal{J}o(2)=1$ from the commutativity of the multiplication,
  \item $\dim \mathcal{J}a\mathcal{J}o(3)=2$  and more precisely $\mathcal{J}a\mathcal{J}o(3)=\K\{x_1(x_2x_3),x_2(x_3x_1)\}$ 
\end{enumerate}
The vector space $\mathcal{J}a\mathcal{J}o(4)$ is generated by the products $(x_i(x_j(x_hx_k))$ and $(x_ix_j)(x_kx_l)$. From the Jacobi-Jordan axiom, we have
$$
\begin{array}{ll}
 (x_ix_j)(x_kx_l)     & =-x_i(x_j(x_lx_k))-x_j(x_i(x_kx_l))   \\
      &   = -x_k(x_l(x_ix_j))-x_l((x_ix_j)x_k)
\end{array}
$$
then $(x_i(x_j(x_kx_l))$ generates $\mathcal{J}a\mathcal{J}o(4)$. Moreover, the previous identity implies
$$
\begin{array}{l}
 x_1(x_2(x_3x_4))+x_2(x_1(x_3x_4))-x_3(x_4(x_1x_2))-x_4(x_3(x_1x_2))=0,    \\
   x_3(x_2(x_1x_4))+x_2(x_3(x_1x_4))-x_1(x_4(x_3x_2))-x_4(x_1(x_3x_2))=0,  \\
   x_4(x_2(x_3x_1))+x_2(x_4(x_3x_1))-x_3(x_1(x_4x_2))-x_1(x_3(x_4x_2))=0.
\end{array}
$$
This shows that $(x_i(x_j(x_kx_l))$ generates $\mathcal{J}a\mathcal{J}o(4)$ for $i=1,2,3$.
Since
$$(x_i(x_k(x_lx_j))=-(x_i(x_j(x_kx_l))-(x_i(x_l(x_jx_k))$$
we deduce, by adding the $3$ identities, that
$$
\begin{array}{l}
x_1(x_2(x_3x_4))+x_1(x_3(x_2x_4))+x_2(x_1(x_3x_4))
+x_2(x_3(x_1x_4))\\
\qquad +x_3(x_1(x_2x_4))+x_3(x_2(x_1x_4))=0, \\
\end{array}
$$
and it is the only relation on this 6 terms then  $\mathcal{J}a\mathcal{J}o(4)$ is the space
$$\K\{(x_1(x_2(x_3x_4)),(x_1(x_3(x_2x_4)),(x_2(x_1(x_3x_4)),(x_2(x_3(x_1x_4)), (x_3(x_1(x_2x_4))\}$$
and
$$\dim \mathcal{J}a\mathcal{J}o(4)=5.$$
Since the product is commutative, $\mathcal{J}a\mathcal{J}o(5)$ is generated by the products represented by
$$\bu(\bu(\bu(\bu\bu))), \ \ \bu((\bu\bu)(\bu\bu)), \  (\bu\bu)((\bu\bu)\bu).$$
The Jordan-Jacobi condition gives (we use now the letters $a,b,c, \cdots$ in place of $x_1,x_2,x_3,\cdots$ to shorten the formulae):
$$\begin{array}{l}
  (ab)((cd)e)+a(b((cd)e)+b(a((cd)e))=0    \\
   a((bc)(de))+(bc)((de)a)+(de)(a(bc))=0  
\end{array}
$$
this shows that $\mathcal{J}a\mathcal{J}o(5)$ is generated by the products $\bu(\bu(\bu(\bu\bu)))$. Since $\dim \mathcal{J}a\mathcal{J}o(4)=5$, we have $25$ generators of this type. We have
$$(ab)((cd)e)+a(b((cd)e)+b(a((cd)e))=0$$
and also
$$
\begin{array}{ ll}
 (ab)((cd)e)&=-(cd)(e(ab))-e((ab)(cd))   \\
      & = c(d(e(ab)))+d(c(e(ab)))+e(a(b(cd)))+e(b(a(cd)))  
\end{array}$$
that implies
$$a(b((cd)e)+b(a((cd)e))+c(d(e(ab)))+d(c(e(ab)))+e(a(b(cd)))+e(b(a(cd)))=0$$
and by symmetry
$$a(c((bd)e)+c(a((bd)e))+b(d(e(ac)))+d(b(e(ac)))+e(a(c(bd)))+e(c(a(bd)))=0.$$
Since the other relations are consequences of the previous relations, we deduce that
$$\dim \mathcal{J}a\mathcal{J}o(5)=23.$$
At this step, the generating series is
$$g(t)=-t+\frac{t^2}{2}-\frac{t^3}{3}+\frac{5t^4}{4!}+\frac{23}{5!}t^5+o(t^6)$$
Computing the ideal of relations of the dual operad $\mathcal{J}a\mathcal{J}o^{!}$, we see that the associated multiplication is anti-commutative and satisfies
$$x_1(x_2x_3)-x_3(x_1x_2)=0$$
The anti-commutativity implies that 
$$x_1(x_2x_3)+(x_1x_2)x_3=0$$
and a $\mathcal{J}a\mathcal{J}o$-algebra is an anti-associative anti-commutative algebra. Then 
$$\mathcal{J}a\mathcal{J}o^{!}(4)=\{0\}.$$
Its generating series is
$$g_{\mathcal{J}a\mathcal{J}o^{!}}=-t+\frac{t^2}{2}-\frac{t^3}{6}$$
and the inverse series
$$-t + \frac{t^2}{2} - \frac{t^3}{3} + \frac{5 t^4}{24} - \frac{t^5}{12} - \frac{7 t^6}{144} + \frac{
 13 t^7}{72} +O(t^8).$$
 Such series cannot be the generating series of a quadratic operad because the sign of $t^6$ is negative. Consequently,  $\mathcal{J}a\mathcal{J}o^{!}$ is not a Koszul operad and
 \begin{theorem}
 The quadratic operad $\mathcal{J}a\mathcal{J}o$ corresponding to the Jacobi-Jordan algebra is not a Koszul operad.
 \end{theorem}
 We deduce, as in the previous case, that the deformation cohomology is described by the cohomology defined by the minimal model.  Since the Poincaré series of  $\mathcal{J}a\mathcal{J}o$ is $g(t)=-t+\frac{t^2}{2}-\frac{t^3}{3}+\frac{5t^4}{4!}+o(t^5)$, the inverse is $\widetilde{g}(t)=-t+ \frac{t^2}{2}-\frac{t^3}{6}+o(t^5)$. We consider the complex
\[
C^1_{\mathcal{J}a\mathcal{J}o}(A,A) \stackrel{\delta^1}{\longrightarrow} C^2_{\mathcal{J}a\mathcal{J}o}(A,A) 
\stackrel{\delta^2}{\longrightarrow}
C^3_{\mathcal{J}a\mathcal{J}o}(A,A)  \stackrel{\delta^3}{\longrightarrow}
C^4_{\mathcal{J}a\mathcal{J}o}(A,A)  \stackrel{\delta^4}{\longrightarrow} \cdots
\]
with
 
-- $C^1_{\mathcal{J}a\mathcal{J}o}(A,A) = Hom(A,A)$

-- $C^2_{\mathcal{J}a\mathcal{J}o}(A,A) = Sym(A^{\otimes^2},A)$

-- $C^3_{\mathcal{J}a\mathcal{J}o}(A,A) = Sym(A^{\otimes^3},A)$, and

-- $C^4_{\mathcal{J}a\mathcal{J}o}(A,A) = Sym(A^{\otimes^5},A) \oplus Sym(A^{\otimes^5},A)  \oplus
 Sym(A^{\otimes^5},A) \oplus Sym(A^{\otimes^5},A)$

\noindent where  $Sym(A^{\otimes^n},A)$ denotes the space of symmetric $n$-linear maps, that is invariant by the action of the symmetric group $\Sigma_n$. The differentials $\delta^i$ of this cohomological complex coincide with the standard cohomology in degree $1$ and $2$:
$$\delta^1(f)(x_1,x_2)=f(x_1x_2)-x_1f(x_2)-f(x_1)x_2$$
and
$$\delta^2 \varphi (x_1,x_2,x_3)= x_1\varphi (x_2,x_3)+\varphi (x_1,x_2x_3)+x_2\varphi (x_3,x_1)+\varphi (x_2,x_3x_1)+x_3\varphi (x_1,x_2)+\varphi (x_3,x_1x_2)$$
where $\varphi$ is a symmetric bilinear form, that is invariant by $\Sigma_2$. We have
$$\delta^2 \varphi (x_1,x_2,x_3)=\sum_{\sigma \in \Sigma_3}\delta^2 _{AA}\varphi (x_{\sigma(1)},x_{\sigma(2)},x_{\sigma(3)})$$
where $\delta^2 _{AA}$ is the differential associated with the cohomological complex of the \ant algebras. 
Then, in degree $4$, we will have
\[
\delta^3 (\psi) = (\sum_{\sigma \in \Sigma_4}\delta_{1,AA}^3(\psi),\delta_{2,AA}^3(\psi),\delta_{3,AA}^3(\psi),\delta_{4,AA}^3(\psi)),
\]
where $\psi$ is invariant by $\Sigma_4$ and $\delta_{i,AA}^3$ corresponding to the \ant case.

Acknowledgment: I thank Liu Xin for useful comments on the first version.

\end{document}